\documentclass[graybox]{svmult}

\usepackage{type1cm}        
\usepackage{makeidx}         
\usepackage{graphicx}        
\usepackage[tight,TABTOPCAP]{subfigure}
\usepackage{multicol}        
\usepackage[bottom]{footmisc}

\usepackage{newtxtext}       %
\usepackage[varvw]{newtxmath}       
\usepackage{tikz,siunitx}

\newcommand{\RN}[1]{%
  \textup{\uppercase\expandafter{\romannumeral#1}}%
}

\newcommand\Rey{\mbox{\textit{Re}}}  

\newcommand{\dif}{{\mathrm{d}}}

\providecommand\bnabla{\boldsymbol{\nabla}}
\providecommand\bcdot{\boldsymbol{\cdot}}
\newcommand{\lapl}[1]{\nabla^2 #1}
\newcommand{\grad}[1]{{\bnabla\! #1}}
\renewcommand{\div}[1]{{\bnabla\! \bcdot #1}}

\newcommand{\qvec}{\boldsymbol{q}}
\newcommand{\Qvec}{\boldsymbol{Q}}

\newcommand{\uvec}{\boldsymbol{u}}
\newcommand{\Uvec}{\boldsymbol{U}}

\newcommand{\xvec}{\boldsymbol{x}}

\newcommand{\uadj}{\uvec^{\dagger}}
\newcommand{\padj}{p^{\dagger}}

\newcommand{\qhvec}{\boldsymbol{\hat{q}}}

\newcommand{\uhvec}{\boldsymbol{\hat{u}}}
\newcommand{\qhadj}{\qhvec^{\dagger}}
\newcommand{\uhadj}{\uhvec^{\dagger}}

\DeclareMathAlphabet\mathsfbi{OT1}{cmss}{m}{sl}

\newcommand\shifthat[2]{%
  \stackengine{\Sstackgap}{$#2$}{\(\hspace{#1}\hat{}\)}{O}{l}{F}{T}{S}}
\newcommand\newhat[1]{%
\if A#1\shifthat{5.2pt}{#1}\else
\if B#1\shifthat{4.8pt}{#1}\else
\if x#1\shifthat{3.6pt}{#1}\else
\shifthat{3.4pt}{#1}
\fi
\fi
\fi
}

\newcommand\No[1][.13ex]{%
  \setbox0=\hbox{\scalebox{.7}{o}}%
  \setbox2=\hbox{n}%
  n\kern-.05em\stackengine{\dimexpr\ht0-\ht2+#1}{\belowbaseline[-\ht2]{\copy0}}%
    {\rule[-.13ex]{.7\wd0}{.13ex}}%
    {U}{c}{F}{F}{L}%
}

\makeindex             

\begin{document}

\title*{Adaptive mesh refinement for global stability analysis of transitional flows}
\titlerunning{AMR for global stability analysis of transitional flows}
\author{Daniele Massaro, Valerio Lupi, Adam Peplinski and Philipp Schlatter }
\institute{D. Massaro \and V. Lupi \and A. Peplinski \and P. Schlatter \at SimEx/FLOW, Engineering Mechanics, KTH Royal Institute of Technology, Stockholm, Sweden. P. Schlatter \at Institute of Fluid Mechanics (LSTM), Friedrich--Alexander--Universität Erlangen--Nürnberg,  Erlangen, Germany. \\  \email{\url{dmassaro@kth.se}}
}
%
%

\maketitle

\abstract{In this work, we introduce the novel application of the adaptive mesh refinement (AMR) technique in the global stability analysis of incompressible flows. The design of an accurate mesh for transitional flows is crucial. Indeed, an inadequate resolution might introduce numerical noise that triggers premature transition. With AMR, we enable the design of three different and independent meshes for the non-linear base flow, the linear direct and adjoint solutions. Each of those is designed to reduce the truncation and quadrature errors for its respective solution, which are measured via the spectral error indicator. We provide details about the workflow and the refining procedure. The numerical framework is validated for the two-dimensional flow past a circular cylinder, computing a portion of the spectrum for the linearised direct and adjoint Navier--Stokes operators. }

\section{Introduction}
\label{sec:dmas_intro}
\begin{figure}
\centering
\includegraphics[trim=8cm 5cm 10cm 3cm,clip,width=0.9\textwidth]{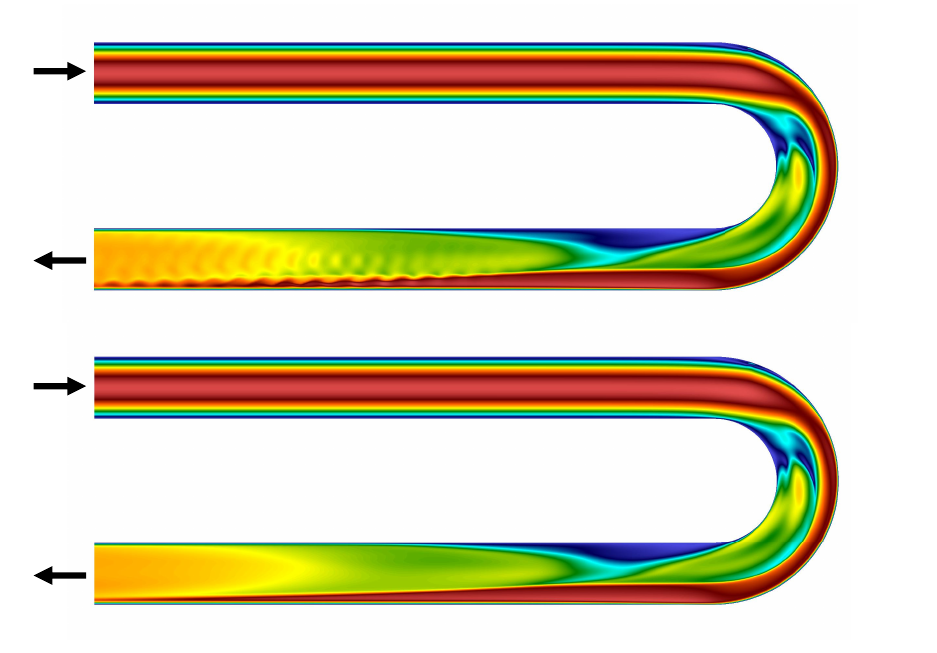}
\caption{Non-linear simulations of three-dimensional spatially developing $180^{\circ}$-bend pipe flow used to extract the base flow. Velocity magnitude in the $x$-$z$ symmetry plane for (\textit{top}) $N = 5$ and (\textit{bottom}) $N = 7$ at $Re_b=2500$ (based on the diameter, the bulk velocity and the kinematic viscosity of the fluid), after 500 convective time units (based on the diameter and the bulk velocity). The black arrows indicate the inflow and outflow directions. This flow configuration is found to be stable up to $Re_b \approx 2550$, however, the low resolution triggers the flow in the outlet at lower Reynolds number, making the base flow extraction diverge. }
\label{fig:num_stab}
\end{figure}
For parallel flows, such as channel and straight pipe flows, classical local stability theory can be applied \cite{Stability}, resulting in one-dimensional problems. Nonetheless, for more complex cases, where the flow is inhomogeneous along all spatial directions, three-dimensional stability analysis, often termed TriGlobal \cite{Theofilis}, needs to be performed. This implies a much higher computational cost and larger memory requirements compared to local analysis. This kind of problem is often tackled using matrix-free, time-stepping algorithms \cite{tuckerman2000} because of the size of the matrices involved. \newline
The dissipation introduced by the discretisation method is a key aspect for accurately capturing the instabilities. Stability analysis is indeed one of the fields where high-order, low-dissipative spatial discretisation methods have been largely applied. Nevertheless, a strong dependency on the polynomial order $N$, \emph{i.e.}\ on the spatial resolution, has been observed even with these methods \cite{peplinski2014}. Thus, the design of the mesh is of paramount importance for stability analysis. As an example, considering the flow through a three-dimensional spatially developing $180^{\circ}$-bend pipe flow, the non-linear direct numerical simulations show numerical instabilities for a low polynomial order (Fig.~\ref{fig:num_stab}). However, generally, it is not known a priori in which regions of the flow a higher spatial resolution is required to accurately capture the instability. For this reason, an adaptive mesh refinement (AMR) technique can be employed. This approach represents an additional tool that allows to refine the computational grid only where needed. In previous works, it has been applied to turbulent flows, see for example Refs.\ \cite{massaro22TSFP,massaro2023step5000,offermans2022}, but there is no study in the literature that employs this technique for stability analysis. In this case, the purpose is not (only) to reduce the computational cost but to provide the proper resolution to robustly capture the instability.\newline
The objective of the present paper is to describe the use of the adaptive mesh refinement technique for global stability analysis. The numerical framework and the steps of the workflow are first described. The procedure is then validated for the case of the two-dimensional flow around a circular cylinder.

\section{Numerical framework}
\label{sec:dmas_mot}

Direct numerical simulations (DNS) are conducted using the open-source code Nek5000 \cite{nek5000-web-page}, which has demonstrated excellent parallel efficiency \cite{Offermans16} and minimal numerical dissipation, crucial attributes for conducting extensive, high-fidelity simulations at a large scale. Regarding the spatial discretisation, the code uses the spectral element method (SEM) \cite{patera1984468}, combining local spectral accuracy, \emph{i.e.}\ nearly exponential error convergence, and isoparametric geometrical flexibility. The computational domain is decomposed into a set of non-overlapping subdomains (elements), and each element is treated as a spectral domain. The $\mathbb{P}_N-\mathbb{P}_{N-2}$  formulation is adopted, \emph{i.e.}\ the functional spaces for the velocity and pressure are spanned by the Lagrangian interpolants integrated over Gauss--Lobatto--Legendre points (GLL) and Gauss--Legendre (GL) points, respectively. The polynomial order is set to $N=7$, corresponding to $N+1$ GLL and $N-1$ GL points. Regarding the time integration, a third-order implicit backward differentiation (BDF) is used, with an extrapolation scheme of order three for the convective term. The advection term is also  de-aliased  with a factor of 3/2 \cite{malm2013}. As the minimum local grid spacing diminishes when an element is repeatedly refined, a variable time step is used in the current simulations. A Courant--Friedrichs--Lewy (CFL) number smaller than $0.3$ is always guaranteed. The tolerances for the velocity and pressure fields are chosen in a conservative way, both below $10^{-8}$. In Nek5000, with the $\mathbb{P}_N-\mathbb{P}_{N-2}$ formulation,  both tolerances are related to the residual in the linear solve, not to the accuracy of the solution, and the pressure tolerance is equal to the desired error in the divergence \cite{nek5000-web-page}. Within this code, our group has implemented \cite{massaro22ico,offermans2019} and extensively used \cite{massaro22TSFP,massaro2023step1000} the adaptive mesh refinement (AMR) technique. The backbone of the AMR implementation in Nek5000 is explained now.  

The first ingredient needed by AMR is an adequate mesh strategy adaptation. Three main possibilities for SEM are available in the literature. The $r$-refinement involves adjusting the boundaries of elements to achieve a desired element size while keeping the number of degrees of freedom (DOFs) constant. This approach is relatively inflexible due to the fixed number of DOFs, which limits the extent of mesh modifications that can be made. The  $p$-refinement involves raising/lowering the polynomial order in specific regions, although Nek5000 does not support the use of elements with varying polynomial orders. Eventually, the $h$-refinement entails subdividing elements, locally refining and coarsening the mesh. This approach proves to be the most flexible form of refinement, especially when there is a substantial increase in the number of elements. Thus the isotropic $h$-refinement was chosen for our applications. As an example, in three-dimensional turbulent external flow around a cylinder with a local discontinuity, the number of elements rises from $4,912$ to $563,652$ spectral elements \cite{massaro2023step5000}. The hanging nodes are not considered as real degrees of freedom, with the conforming-space/non-conforming-mesh approach. The non-conforming interfaces are treated with the parent-to-children interpolating operator \cite{kruse1997}, also assessing that non-conforming interfaces do not affect the solution or introduce any instabilities \cite{massaro22ico}. For managing the hierarchical structure of the mesh grid and executing parallel partitioning, we depend on external libraries such as p4est \cite{p4est} and parMetis \cite{parmetis}, or alternatively, parRSB.

The second, and last, piece of the puzzle is a proper estimation of the error. Various measures of the error exist, using either local or global goal-oriented errors. In previous works, we provide a detailed comparison between local error indicators and dual-weighted goal-oriented error estimator, in both steady and unsteady flows \cite{offermans2020,offermans2022}. In the current work, we rely solely on the spectral error indicator (SEI), introduced below.

\subsection{Spectral error indicator}

In this context, we use the spectral error indicator (SEI) by Mavriplis (1989) \cite{mavriplis1989} as a measure of the mesh quality. The SEI is a cost-effective and localised method used to assess the truncation and quadrature errors in the solution field. It falls under the category of indicators, as it offers insights into the current solution's accuracy. In contrast, other estimators involve a rigorous mathematical assessment of error bounds by utilising an alternate solution, such as the adjoint (dual) solution \cite{offermans2022}. For the sake of simplicity, let us consider a one-dimensional problem, where $u$ is the exact solution to a system of one-dimensional partial differential equations, and $u_N$ is an approximate spectral element solution with polynomial order $N$. We expand $u(x)$ on a reference element in terms of the Legendre polynomials:
\begin{equation}
u(x) = \sum_{k=0}^{\infty} \hat{u}_k L_k(x)
\label{eq_1}
\end{equation}
where $\hat{u}_k$ are the associated spectral coefficients and $L_k(x)$ is the Legendre polynomial of order $k$. The estimated error $\epsilon = \| u-u_N \|_{L^2}$ results in:
\begin{equation}
\epsilon = \Biggl( \int_{N}^{\infty} \frac{\hat{u}(k)^2}{\frac{2k+1}{2}} dk + \frac{\hat{u}^2_N}{\frac{2N+1}{2}} \Biggr)^{\frac{1}{2}}
\label{eq_2}
\end{equation}
where we assume an exponential decay for the spectral coefficients of the form $\hat{u}(k) \approx c \hspace{0.5mm} \text{exp}(-\hat{\sigma} k)$. The parameters $c$ and $\hat{\sigma}$ are obtained interpolating in a linear least-squares sense the $\text{log}(\hat{u}_k)$ w.r.t.\ $k$ for $k\leq N$. In a multi-dimensional problem, the maximum error among each component is considered, providing a single spectral indication per element. As discussed in previous works \cite{offermans2020, offermans2022}, this indicator is well suited to track flow features with high gradients, such as shear layers or fluid-wall interaction regions. Contrary to other works \cite{huet2023}, we do not aim to track instantaneous flow features, but rather to converge to a statistically stationary mesh. Thus, the time-average of $\epsilon$ for a given interval $T$ is calculated. As Nek5000 is mainly written in Fortran 77, dynamic memory allocation is not allowed. The maximum number of elements for a given simulation is specified in the compiling phase. Due to this constraint, the criterion used to define the number of elements to refine at each round is based on a given percentage of the total. As an example, in the current simulations, at each round, the number of elements is increased by 20\% and the threshold on the SEI (to decide if to refine or not a certain element) is based on such restriction. Further details about the refinement strategy and sampling interval are provided below.

Eventually, it is worth noting that the current study required few, but significant, changes in the code. Particularly, the AMR has been adapted for the linear solver in Nek5000. As previously done for the non-linear part, the enforcement of the continuity constraint at the boundaries of the elements required to be modified. The basis functions need to be in the Sobolev space $H^1$, but, at the non-conforming interfaces, the presence of a larger number of degrees of freedom on one side does not allow for strict enforcement of the continuity constraint. Because of this, it was required to interpolate the solution from the side with the lower resolution (parent element), onto the side with the higher resolution (children elements). Thus, the continuity is again enforced by extrapolating the solution from the parent element onto the children elements. A further change was related to the measure of the error since the SEI is estimated by measuring the truncation and quadrature error on the velocity perturbation field rather than on the velocity field itself. Eventually, the AMR compatibility with the stability tools provided by the KTH framework \cite{KTHframework} was also required, \emph{e.g.}\ the support for non-conforming mesh in the matrix-free time-stepper approach for the calculation of the eigenvalues \cite{bagheri_matrix}.

\section{Global stability analysis}
\label{sec:dmas_stab}
Global stability analysis investigates the evolution of infinitesimal perturbations on top of a reference solution called base flow. In the case under investigation, the base flow $\Qvec = (\Uvec, P)^T$ is an exact solution of the steady Navier--Stokes equations
\begin{equation}
\left\{
\begin{array}{ll}
\displaystyle \left( \Uvec  \bcdot \bnabla \right) \Uvec  = - \grad{P} + \dfrac{1}{\Rey} \lapl{\Uvec}, \\[8pt]
\div{\Uvec} = 0,
\end{array}
\right.
\label{eq:ns}
\end{equation}
which are made dimensionless by scaling with a reference velocity $U$, a reference length $L$ and the constant density of the fluid $\rho$. The Reynolds number is thus defined as $\Rey = U L/\nu$, where $\nu$ is the kinematic viscosity.\newline
The first step of the global stability analysis consists of computing the steady base flow. If it is stable, it can be easily computed by integrating in time the Navier--Stokes equations. However, when an unstable steady state is sought, convergence to the base flow is obtained by providing additional stabilisation, \emph{e.g.}\ through selective frequency damping (SFD) proposed by Ref.\ \cite{akerviketal2006}.\newline
Since infinitesimal disturbances $\qvec' = (\uvec', p')^T$ are considered, their dynamics can be described by linearising the Navier--Stokes equations about the base state, reading
\begin{equation}
\left\{
\begin{array}{ll}
    \dfrac{\partial \uvec'}{\partial t} + (\uvec' \bcdot \bnabla)\ \Uvec + (\Uvec \bcdot \bnabla)\ \uvec'=-\grad{p'} + \dfrac{1}{\Rey} \lapl{\uvec'}, \\[8pt]
    \nabla \cdot \uvec'=0,
\end{array}
\right.
\label{eq:LNSE}
\end{equation}
The solution of the direct linear problem \eqref{eq:LNSE} converges to the least stable eigendirection in the asymptotic limit \cite{Strogatz}. Introducing the normal-mode ansatz
\begin{equation}
    \uvec'(\xvec,t) = \uhvec(\xvec) \:\! e^{\lambda \:\! t}, \quad p'(\xvec,t) = \hat{p}(\xvec) \:\! e^{\lambda \:\! t}
    \label{eq:triglobal_ansatz}
\end{equation}
into equations \eqref{eq:LNSE} leads to the generalised eigenvalue problem
\begin{equation}
    \lambda \:\! \mathcal{R} \:\! \qhvec = \mathcal{L} \:\! \qhvec,
    \label{eq:gen_eigval_problem}
\end{equation}
where $\mathcal{R}$ is a restriction operator  (identity for the velocity, zero for the pressure), $\mathcal{L}$ is the linearised Navier--Stokes operator, $\qhvec = (\uhvec, \hat{p})^T$ is a complex-valued eigenfunction and $\lambda = \sigma + \imag \, \omega \in \mathbb{C}$ is the associated eigenvalue, with $\sigma$ being the growth rate and $\omega$ the angular frequency. The base flow is unstable to infinitesimal disturbances if there exists at least one eigenvalue with $\sigma > 0$. The critical Reynolds number $\Rey_{cr}$ is defined as the value for which the base flow is marginally stable, \emph{i.e.}\ $\sigma = 0$. \newline
Adjoint equations are used in stability analysis to provide the sensitivity of an objective to different inputs \cite{luchini2014}. The adjoint global eigenmodes are solutions to the generalised eigenvalue problem that can be derived from the adjoint linearised Navier--Stokes equations
\begin{equation}
\left\{
\begin{array}{ll}
-\dfrac{\partial \uadj}{\partial t} - \left(\Uvec \bcdot \bnabla \right) \uadj + \left(\bnabla \Uvec \right)^\mathrm{T} \uadj  =  \grad{\padj} + \dfrac{1}{\Rey_b} \lapl{\uadj}, 
\\[8pt]
-\div{\uadj} = 0,  
\end{array}
\right.
\label{eq:adjointNS}
\end{equation}
with $\uadj$ and $\padj$ being the adjoint perturbation velocity and pressure, respectively. The receptivity of a direct mode to initial conditions and momentum forcing is given by the velocity field of the corresponding adjoint eigenmode $\qhadj = (\uhadj, \hat{p}^{\dagger})^T$, whereas $\hat{p}^{\dagger}$ represents its sensitivity to mass sources \cite{giannetti2007}. As for the direct problem, the solution of the system \eqref{eq:adjointNS} approaches the least stable adjoint eigenmode for $t \rightarrow \infty$. Note also that the adjoint spectrum is the complex conjugate of the direct one \cite{robinson_2020}. 

\section{Methodology \& Results}

Generally, the global stability analysis of a transitional flow requires the solution of three different sets of equations, as commented above. Each of these steps requires a mesh that can be significantly different from the others. Thus, we aim to introduce the usage of AMR in the study of transitional flows, where the mesh for each solution is designed to minimise the numerical error on that specific solution field. The general procedure is described below. 

\begin{figure}
\centering
\includegraphics[trim=0cm 0cm 1.2cm 2cm,clip,width=\textwidth]{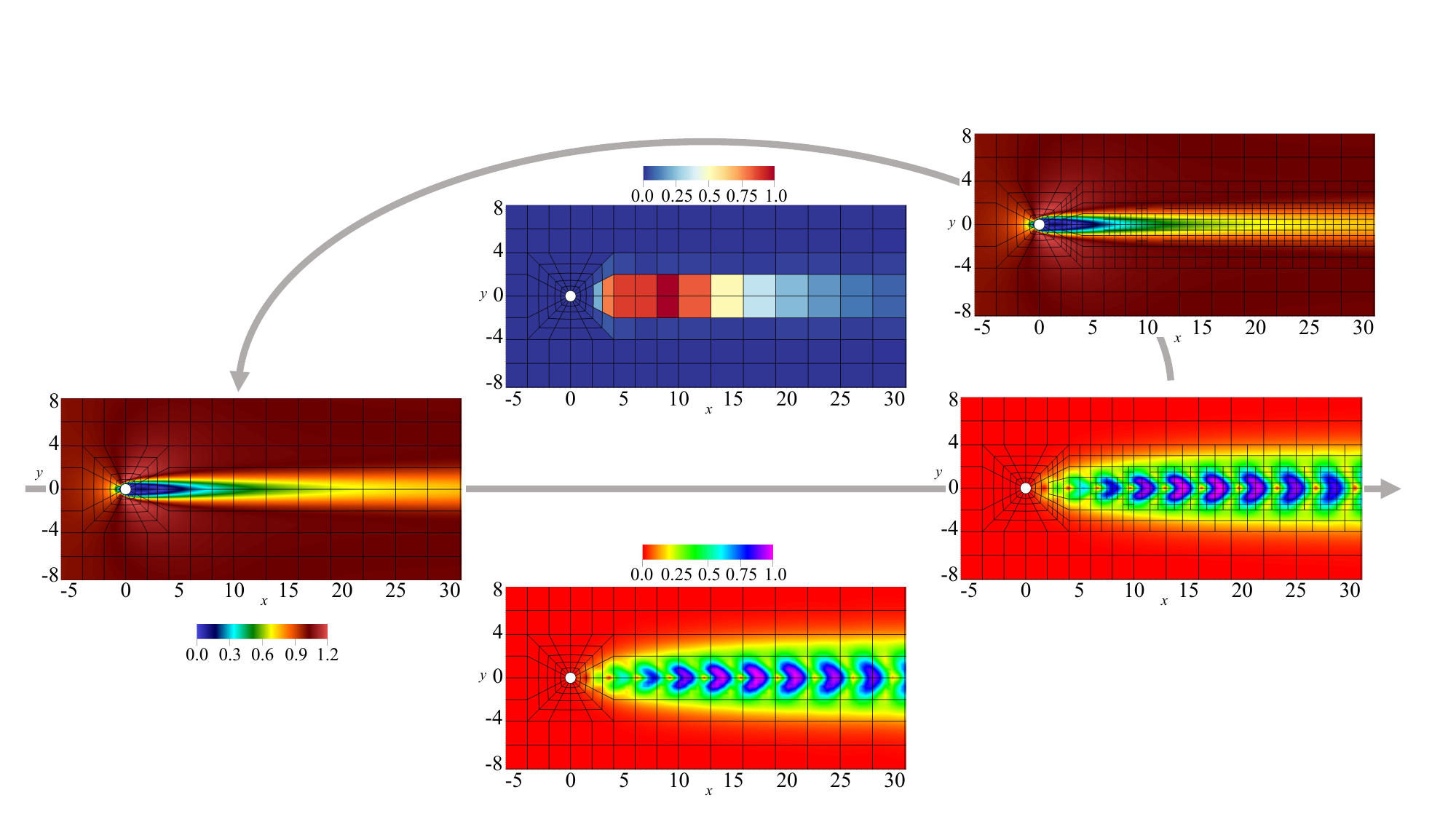}
\caption{Workflow of the AMR application for the global stability analysis of the flow past a circular cylinder at $Re_D=50$. The base flow extracted via SFD is spectrally interpolated onto a conforming mesh (Mesh A2). The solution, after being initialised with spatially uncorrelated noise, converges to the most unstable eigenvector. The SEI is measured on the perturbation velocity field. The mesh is refined according to the SEI. The procedure is repeated several times. The base flow (which is frozen when integrating the linearised Navier--Stokes equations) is also interpolated on the refined mesh since Nek5000 does not allow handling different meshes. The spectral error indicator and the linear perturbation fields are normalised w.r.t.\ their maximum value. }
\label{fig:flow}
\end{figure}

\begin{enumerate}
    \item Non-linear direct numerical simulations
    \begin{enumerate}
         \item[a)] A set of initial simulations is carried out to identify, approximately, the critical Reynolds number $Re_{cr}$.         
     \end{enumerate}
    \item Base flow calculation
    \begin{enumerate}
         \item[a)] The incompressible (non-linear) Navier--Stokes equations are numerically integrated for a Reynolds number slightly above $Re_{cr}$. Selective frequency damping \cite{akerviketal2006} is used to stabilise the flow and extract the base flow.
         \item[b)] After a short initial transient, the SEI is collected on the velocity field ($\Uvec$-variable). The collection lasts for a given (and case-dependent) time interval $T$. Then, the SEI is time-averaged, and the mesh is refined accordingly.
         \item[c)] Multiple levels of refinement are performed until an adequate resolution is obtained. 
         \item[d)] At this stage, the mesh is frozen (Mesh A), and the simulation ends when the SFD tolerance converges below a given threshold.
     \end{enumerate}
    \item Direct/Dual linear solution
    \begin{enumerate}  
         \item[a)] The linear simulation is initialised with the mesh (Mesh A) that was designed for the base flow solution. Alternatively, Mesh A is spectrally interpolated on a new mesh (Mesh A2), which becomes the initial spectral grid for the linear simulation. The interpolation is necessary whether we want to start the linear simulation on a fresh spectral grid. Particularly, we perform a spectral interpolation,  with an accuracy that corresponds to the adopted polynomial order ($p=7$) \cite{azadinterp}. The effects of the mesh reduction and the base flow interpolation have to be carefully assessed to make sure the instability mechanism is not altered by them.    
         \item[b)] The linearised direct/dual Navier--Stokes equations are initialised with noise uncorrelated in space which has a non-zero projection on the wanted modes and a frozen base flow (previously extracted via SFD).
         \item[c)] The disturbance converges to the (globally) most unstable eigenmode.        
         \item[d)] After a short (case-dependent) initial transient, the collection of the SEI begins. The choice of the transient $T_d$ is crucial, especially for the linear simulations. Indeed, on the one hand, a short $T_d$ can lead to the over-refinement of some regions, not necessary for the converged unstable eigenmode. On the other hand, when $T_d$ is too long, the still coarse mesh can trigger the flow transition due to the numerical error. A conservative choice of $T_d$ is preferable, as our framework also allows to coarse the grid, whether this is necessary. The SEI is collected on the perturbation velocity field ($\uvec'$-variable) and time-averaged for a reasonable time interval $T$. The mesh is refined accordingly. Note that the perturbation field is exponentially growing in time, thus a rescaling might be required. However, such a scaling factor does not affect the outcomes of the refinement process. 
         \item[e)] The procedure is repeated several times until an adequate resolution is obtained. Standard convergence analysis can be carried out by looking at local (velocity probes) and global (total perturbation energy) quantities. When the convergence is assessed, the mesh is frozen (Mesh B).
         \item[f)] Mesh B is used for the eigenvalues (and corresponding eigenvectors) calculation of the direct/dual linearised Navier--Stokes operator (see Section \ref{sec:dmas_stab}) via the Implicitly Restarted Arnoldi Method (IRAM). The stability tools by Ref.\ \cite{KTHframework} have been adapted to handle non-conforming (and, at this stage, static) meshes. 
         \item[g)] Steps from a) to e) are repeated for the dual linear solution.       
    \end{enumerate}

\end{enumerate}

\subsection{Flow past circular cylinder}

The global stability of the external flow around a circular cylinder with mesh adaptivity is considered. The cylinder has a unit diameter $D$ and is placed at $(0,0)$ in a spectral grid extending from $-15D$ to $35D$ and $-15D$ to $15D$ in $x$ and $y$ directions, respectively. A Dirichlet boundary condition (uniform unitary velocity $U$) is set at the inlet, and periodic boundary conditions are used at $y=\pm 15D$. At the outflow, a natural boundary condition $ (-p \mathbf{I} + \nu \nabla \mathbf{u})\cdot \mathbf{n} = 0 $, where $\mathbf{n}$ is the normal vector, is imposed. The polynomial order is $N=7$, and the initial (coarse) mesh counts approximately $300$ elements. 

\begin{figure}
\centering
    \subfigure[]
    {\includegraphics[trim=9.5cm 2cm 7.5cm 5cm,clip,width=0.495\textwidth]
        {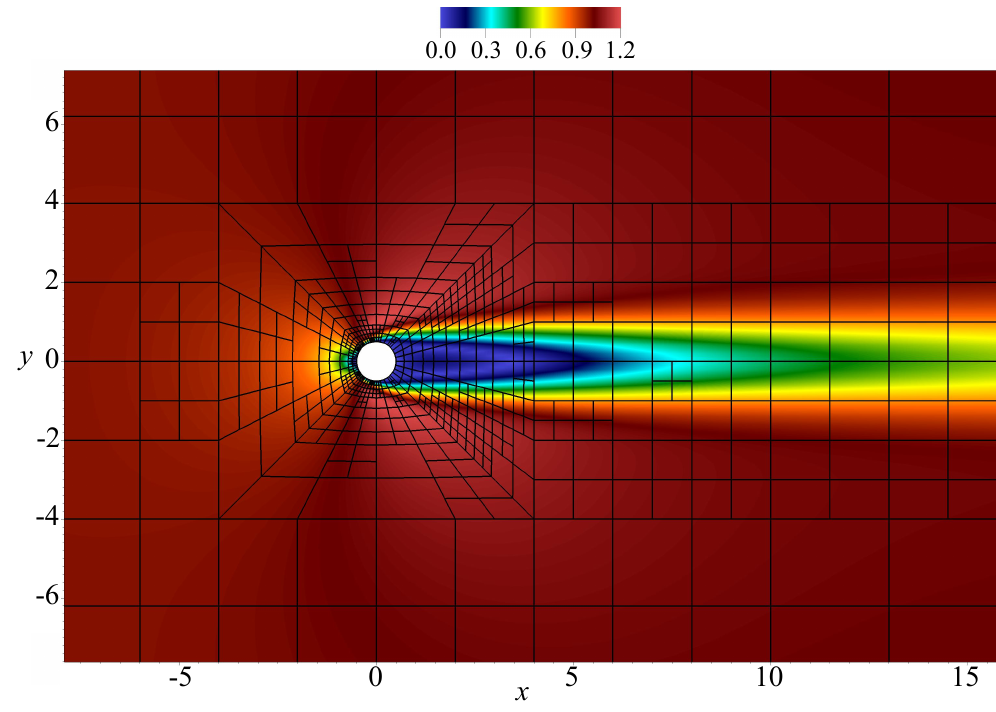}}

    \subfigure[]
    {\includegraphics[trim=9.5cm 2cm 7.5cm 5cm,clip,width=0.495\textwidth]
        {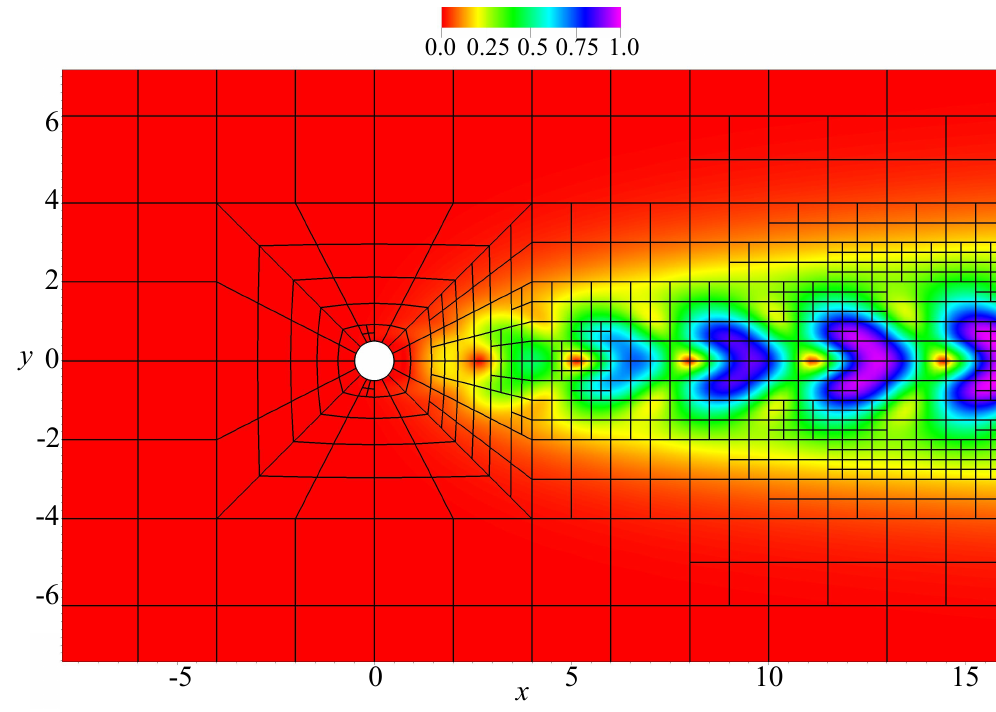}} 
    \subfigure[]
    {\includegraphics[trim=9.5cm 2cm 7.5cm 5cm,clip,width=0.495\textwidth]
        {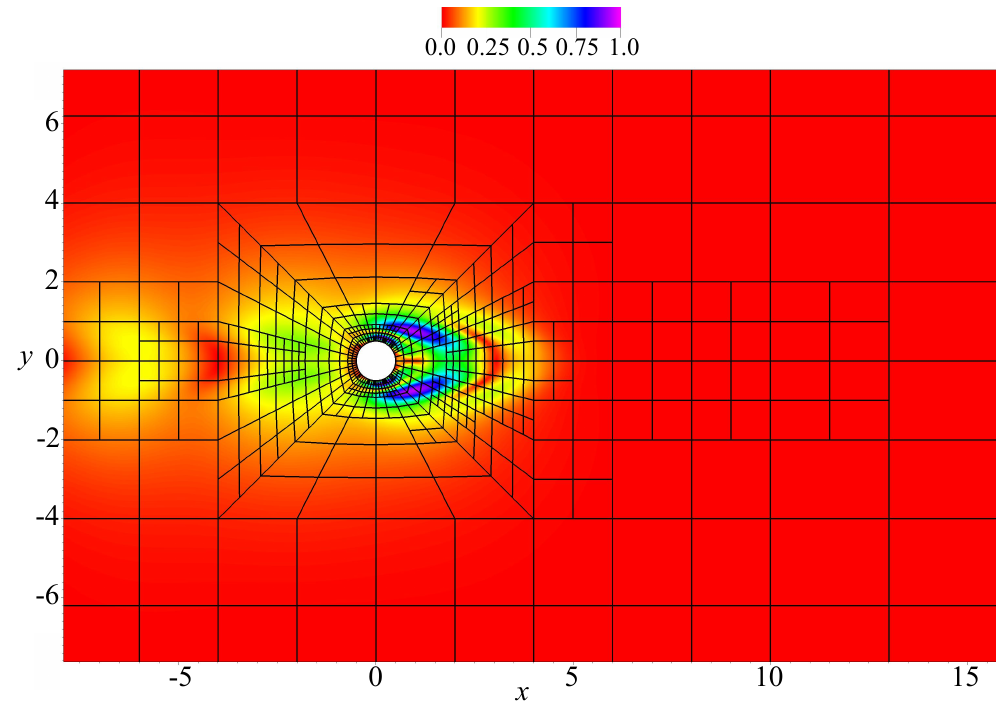}}
    \caption{The base flow (a), the linear direct (b) and dual (c) most unstable eigenmode with a spectral grid designed for each of those different solutions. Only a portion of the domain is shown, and the linear direct and dual solutions are normalised by their maximum. }
    \label{fig:global}
\end{figure}

\subsection{Non-linear base flow calculation}

\begin{figure}
\centering
\includegraphics[trim=2.5cm 1cm 2.5cm 1.5cm,clip,width=\textwidth]{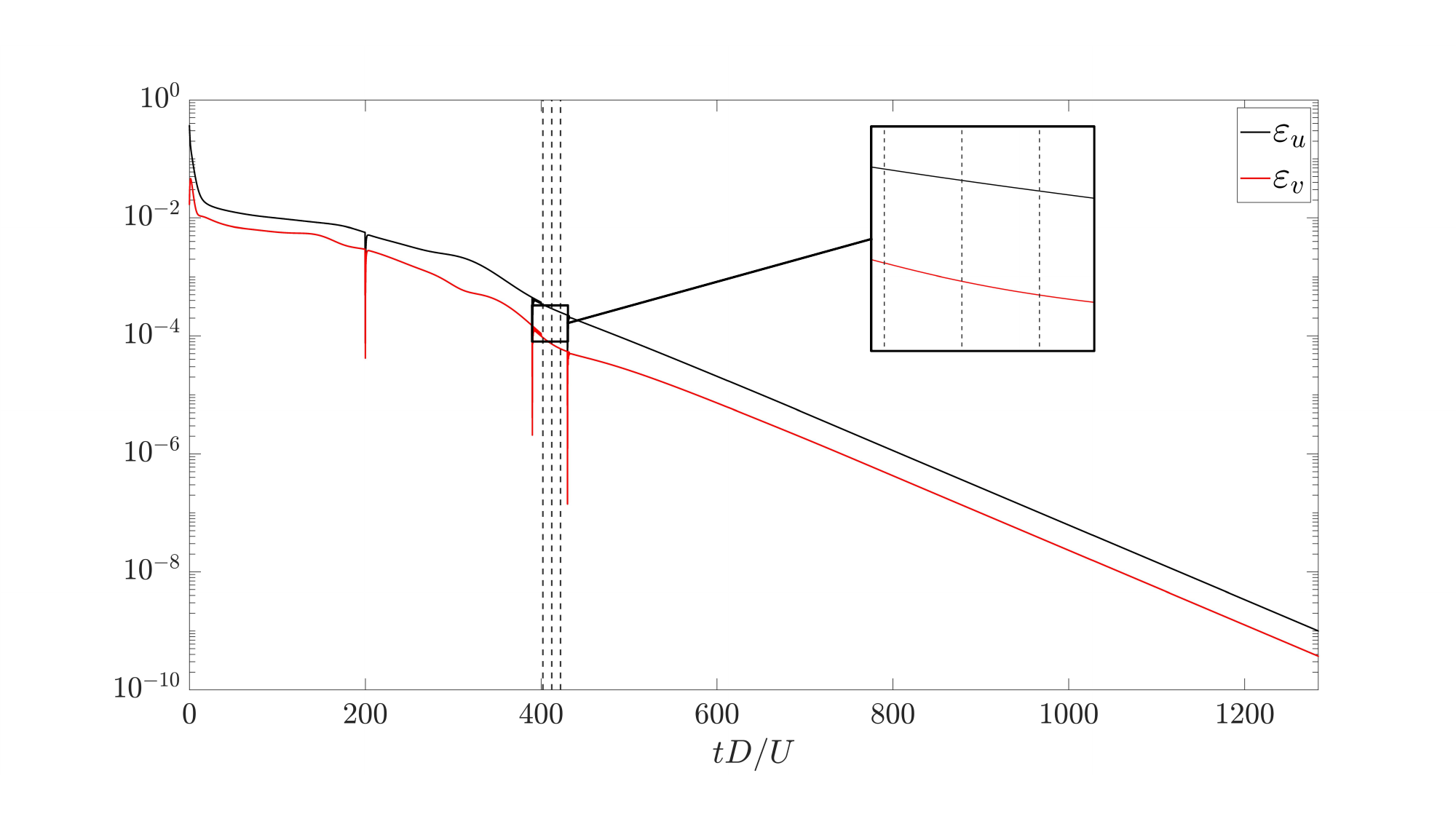}
\caption{The convergence of the selective frequency damping (SFD) tolerance $\varepsilon$ for the two velocity components. The time instants when the refinement is performed are indicated by vertical dashed lines. The vertical oscillations are related to the simulation restarts. }
\label{fig:SFD}
\end{figure}

For the flow past a circular cylinder, the critical Reynolds number is well-known to be $Re_D \approx 47$ \cite{Tritton,Williamson96}, based on the cylinder diameter $D$ and the inflow velocity $U$. Thus, the initial non-linear direct numerical simulations to identify the critical Reynolds number range are not necessary. The base flow is extracted via the selective frequency damping (SFD) technique \cite{akerviketal2006} at a Reynolds number slightly higher ($Re_D=50$) than the critical value. The SFD damps the oscillations of the solution using a temporal low-pass filter by applying to the flow a forcing $\boldsymbol{f} = -\chi(\boldsymbol{u}-\boldsymbol{w})$, where $\boldsymbol{u}$ is the flow solution and $\boldsymbol{w}$ is the temporally low-pass filtered velocity obtained by a differential exponential filter $\boldsymbol{w}_t=(\boldsymbol{u}-\boldsymbol{w})/\Delta$, with $\Delta$ determining the filter width. Here, according to the formulation adopted in the KTH framework \cite{KTHframework}, we choose $\Delta=4.05$ and $\chi=0.5$. The amplitude of the forcing $\varepsilon = \| \boldsymbol{u}-\boldsymbol{w} \|_{L^2(\Omega)}$ is considered as convergence indicator. The simulation is stopped when the tolerance $\varepsilon$ is sufficiently low to not affect the calculation of the eigenvalues, $\varepsilon <10^{-10}$ in our case. The SFD convergence is shown to be robust when the refinement is performed, see figure \ref{fig:SFD}. No jump appears in the convergence of $\varepsilon$ at the instants when the spectral grid is adapted. A few oscillations are visible at the restart among different runs (dashed lines in Fig.\ \ref{fig:SFD}). Nonetheless, these do not appear to be in any way related to the adaptive mesh process.

After the initial transient, we perform multiple levels of refinements. The SEI is not an absolute measure, but rather an indicator, \emph{i.e.}\ the value itself is meaningless. For this reason, here, we report the normalised SEI in Fig.\ \ref{fig:flow}. However, for the interested reader, typical largest values for the SEI are approximately $10^{-6}$ and $10^{-4}$, in two- and three-dimensional cases, respectively (as reported by Refs. Offermans \textit{et al.} (2020, 2023) \cite{offermans2020,offermans2022}). The SEI is measured for each component of the velocity field $\Uvec$. The maximum among the velocity components is taken, time averaging for an interval $T$. In the current case, the time interval $T$ is set as a multiple of the vortex shedding period. The final mesh is shown in figure \ref{fig:global}. The results align with previous findings, when the AMR was applied to steady flows \cite{offermans2020}. Namely, the SEI is sensitive to the perturbations caused by the cylinder and convected by the bulk flow, focusing the refinement around and close to the surface of the cylinder and propagating downstream, in the wake and regions of higher velocity gradients.

\subsection{Linear direct and dual most unstable eigenmodes}

\begin{figure}
\centering
\includegraphics[trim=2.5cm 1cm 2.5cm 1.5cm,clip,width=\textwidth]{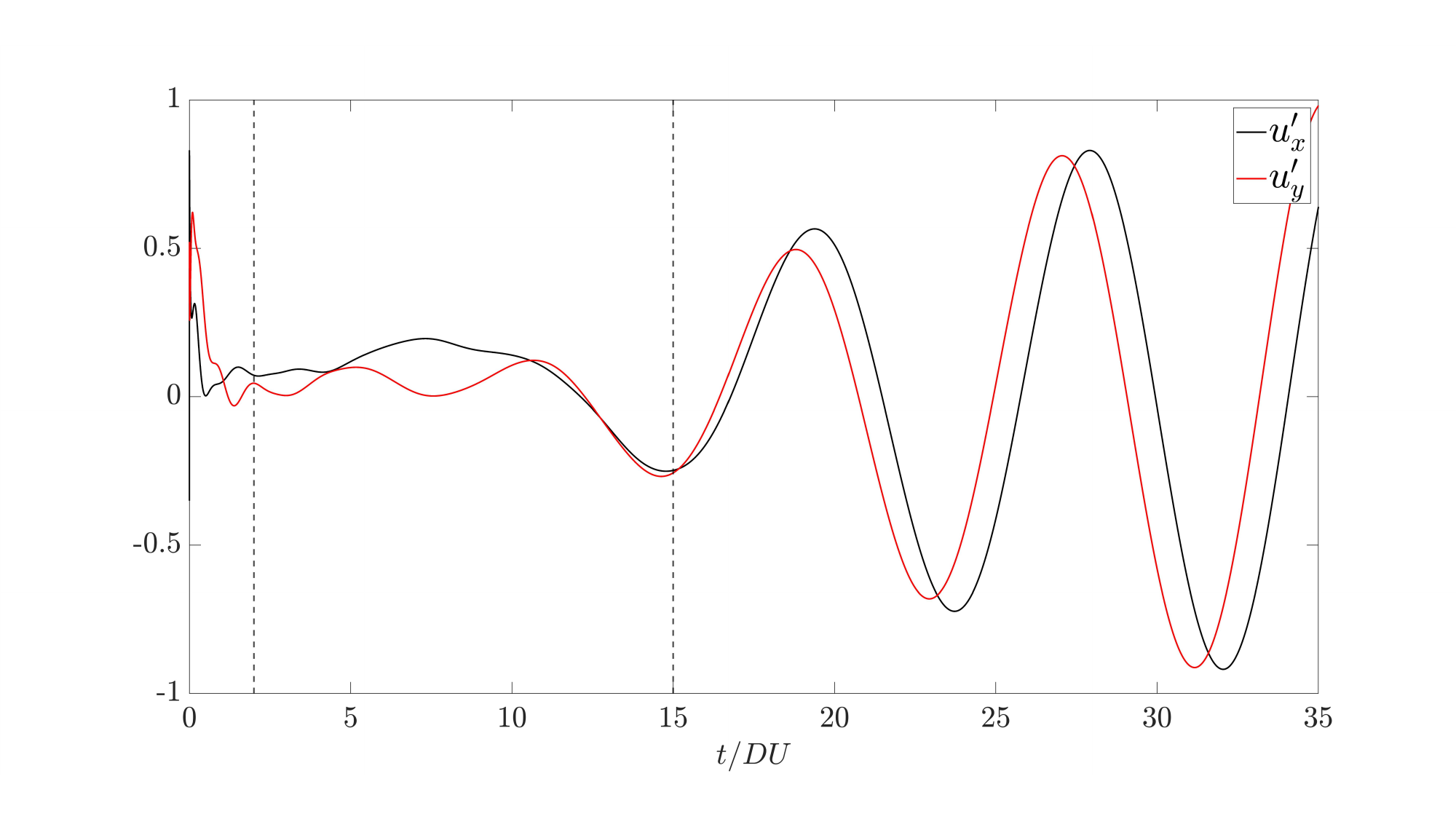}
\caption{Time signals of local velocity probes for the linear direct solution $\uvec'$ in the cylinder wake at ($x/D=6,y/D=1$). The dashed lines indicate two possible choices for the transient time interval $T_d$, before starting the refinement. The solutions are normalised w.r.t.\ their maximum.}
\label{fig:probes}
\end{figure}

The extracted base flow is spectrally interpolated onto a new mesh (Mesh A2) which is used as the initial mesh to design the spectral grids of the linear direct and dual solutions. The Navier--Stokes equations are linearised about the extracted base flow and initialised with spatially uncorrelated disturbance. Marching in time, the disturbances tend to converge to the most unstable eigenmode and after an initial transient $T_d$, the SEI is calculated for the velocity perturbation fields. The error-driven design of the meshes is stopped when the unstable eigenmode is converged to its final spatial shape. Standard mesh convergence analysis are carried out. It is worth noting that the choice of $T_d$ is critical, and varies on a case-to-case basis. In the current simulation, the time signals of local velocity probes in the cylinder wake help define an adequate $T_d$. In figure \ref{fig:probes}, time signals of the two velocity components of the linear direct solution are shown at ($x/D=6,y/D=1$). The first dashed line, corresponding to 2 time units (expressed in terms of $U$ and $D$), represents a conservative choice. The high oscillations due to the development of the initial noise are gone, but the shedding frequency has not been established yet. Generally, it is convenient to start the refinement after such a time interval. The only possible drawback of an early-stage refinement is to build up a mesh that is over-refined in some transient regions. This can be overcome by turning on the coarsening, which was implemented as well. Alternatively, a larger $T_d$ can be considered. As an example, the second dashed line corresponds to 15 time units. At this stage, the perturbation oscillates with a frequency almost identical to the final vortex shedding frequency. The exponential growth is also visible. However, this approach has one main concern: letting the initial disturbance evolve for a longer time on a coarse mesh can introduce some numerical errors that trigger the transition. For this reason, a conservative choice is taken here, setting $T_d=2$. Similarly, rather than local velocity measurements, a global quantity, alike the total perturbation energy 
\begin{equation}
    E_{pert} = \int_V \rho (\uvec')^H \uvec' \dif V,
\end{equation}
could be considered, with the superscript $H$ denoting the conjugate transpose. In the initial transient, this can have very high oscillations, depending on the parameters of the initial condition. When its incipient monotonic growth appears, the mesh refinement can begin.

The final meshes and solutions are shown in figure \ref{fig:global}. The differences among those are remarkable, particularly between the linear and dual eigenmode. Indeed, the non-normality of the Navier--Stokes operator leads to a noteworthy spatial separation between them, requiring different levels of accuracy in different regions. The solution of the direct mode is skew-symmetric w.r.t. the centreline of the wake, and the real and imaginary parts are identical, apart from a given phase shift. The adjoint mode identifies regions of maximum receptivity. These are localised in the near wake of the cylinder, close to the upper and lower sides of the body surface. In contrast to the direct mode, the receptivity decays rapidly both upstream and downstream of the cylinder \cite{giannetti2007}. A preliminary validation is performed by computing the growth rate $\sigma$ for the linear unstable eigenmode with the final non-conforming and a well-resolved conforming mesh. In both cases, the growth rate is $\sigma \approx 0.016$, at $Re_D=50$, with an accuracy $\approx 10^{-6}$. The critical Reynolds number corresponds to the well-known $Re_{cr} \approx 47$. Further validations are carried out by looking at the eigenvalues of both operators in the next section.

\subsection{Spectrum of the linearised operators}

A portion of the spectrum of eigenvalues for both the direct and adjoint linearised Navier--Stokes operator is computed. We reformulate the linearised equations as an initial value problem for the velocity only by exploiting the incompressibility constraint \cite{schlatter_2011} and the eigenpairs are computed with the implicitly restarted Arnoldi method (IRAM), proposed by \cite{Sorensen_1992,ARPACK}. In the Arnoldi method, the eigenpairs are estimated by exploring solutions within a Krylov subspace, specifically with a dimensionality of $m=200$. We calculate 50 eigenpairs for both the primary and associated adjoint problems, employing a residual tolerance of $10^{-6}$ for the eigenvalue calculation. Further details are available in  Refs.\ \cite{bagheri_matrix,KTHframework,peplinski_2015}.

\begin{figure}
\centering
\includegraphics[trim=0cm 0cm 0cm 0cm,clip,width=\textwidth]{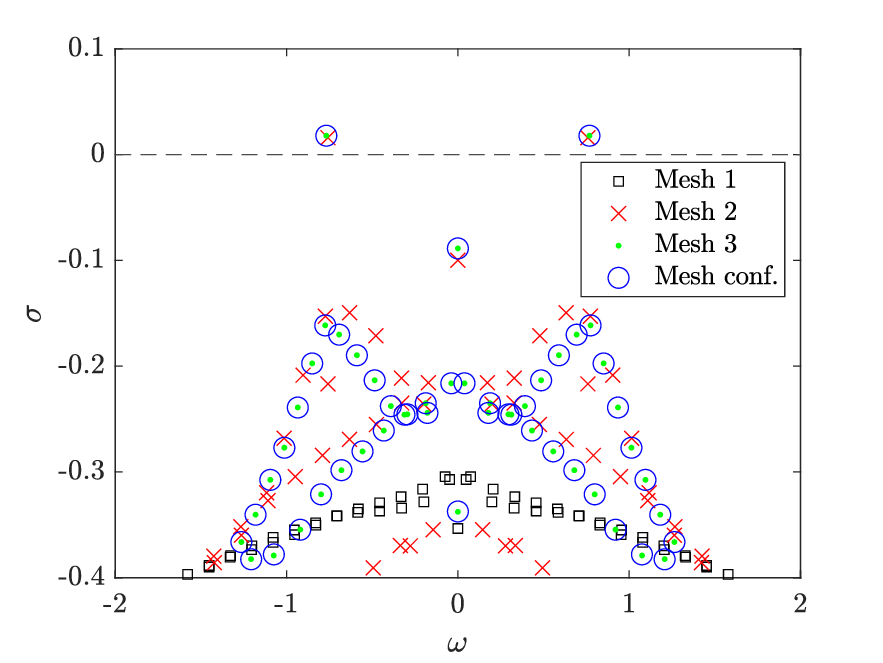}
\caption{Portion of the spectrum of eigenvalues for the linear direct operator computed via the implicitly restarted Arnoldi method. Comparison among three progressively refined meshes with approximately 300, 1200 and 3000 spectral elements, and a well-resolved conforming mesh.}
\label{fig:spec}
\end{figure}

A portion of the spectrum of the linearised (direct) Navier--Stokes operator is shown in figure \ref{fig:spec}. The initial conforming mesh (Mesh 1) with approximately 300 spectral elements ($\sim 20,000$ GLL points) is progressively refined, following the strategy described above. Mesh 2 and 3 contain around 1200 and 3000 elements. Reference data from a well-resolved (5000 elements, $\sim 320,000$ GLL points) conforming mesh are also included. It is clear that the coarsest mesh does not guarantee an adequate resolution. Indeed, not only it does not capture the most unstable eigenmode, but also shows a slight difference between the direct and adjoint eigenvalues, see table \ref{tab:comp}. These latter are supposed to be identical according to linear algebra theory \cite{robinson_2020}. Mesh 2 already captures the most unstable growth rate and angular frequency with an accuracy of $10^{-4}$. Small deviations are visible in the stable eigenvalue pairs. However, these are not crucial for the global stability analysis. Eventually, the results obtained with Mesh 3 are in excellent agreement with the reference conforming grid. The difference for both $\sigma$ and $\omega$ is below $10^{-9}$. Table \ref{tab:comp} confirms the accuracy of the results by also reporting the values for the most unstable dual eigenmode.

\begin{table}
  \begin{center}
\def~{\hphantom{0}}
\begin{tabular}{c|c c c c c }
& $\sigma_{D}$ & $\omega_{D}$ & $\sigma_{A}$ & $\omega_{A}$  \\
\hline
Mesh 1      & $-0.30452389$	& $0.07426587$ & $-0.305515448$ & $0.07389150$  \\
Mesh 2      & $0.01746561$	& $0.76810699$ & $0.017843084$ & $0.76810422$  \\
Mesh 3      & $0.017846528$	& $0.76800714$ & $0.017843042$ & $0.76800439$  \\
Mesh conf.  & $0.017846529$ & $0.76800715$ & $0.017843064$ & $0.76800402$  \\
\end{tabular}
\caption{Growth rate $\sigma$ and angular frequency $\omega$ for the most unstable eigenmode of the direct ($D$) and adjoint ($A$) operator. Comparison among the initial coarse mesh, two non-conforming meshes and a well-resolved conforming mesh used as reference.}
\label{tab:comp}
\end{center}
\end{table}

\section{Conclusions \& Outlooks}
\label{sec:dmas_conclusion}

The current study extends the usage of adaptive mesh refinement to global stability analysis of transitional flows. The numerical framework and the AMR implementation are introduced. We explain the methodology and the crucial aspects, \emph{e.g.}\ the refinement based on different variables for different solutions. Our AMR implementation shows robustness, \emph{i.e}\ the non-conforming interfaces do not introduce any further instabilities in the flow solution, as previously observed for turbulent flows \cite{massaro22ico}. The application for the global stability analysis of a flow past a circular cylinder is described. Three different and independent meshes are designed by minimising the committed truncation and quadrature errors. These latter are measured via the spectral error indicator. Reporting a portion of the spectrum of the linear direct operator (Fig.\ \ref{fig:spec}), we show that a badly resolved mesh can lead to the detection of the wrong physical mechanism. In the worst scenarios, the transition is either anticipated or postponed, as it occurs in the current case. The improvement of the resolution by performing multiple rounds of refinement led to a spectral grid able to capture the unstable eigenmode with an accuracy of $10^{-9}$ and using half of the grid points. The computational saving is expected to be even larger for more complicated or three-dimensional cases. An extension of this work for three-dimensional spatially developing bent pipe flow is described in Ref.\ \cite{massaroStabBent}.

Future development of potential interest is the capability to handle different meshes into Nek5000. This would allow, for example, to have multiple grids for each of the first $n$ eigenmodes in the Arnoldi method, which could be important when multiple unstable eigenvalues exist. Nonetheless, this would require significant modifications to the solver, not only in the AMR implementation but also to the core of the code. It is also worth mentioning the need for an adequate refinement criterion. So far, we rely on standard mesh convergence techniques. However, it would be appealing to have a criterion based on the quality of the estimated eigenvalues which automatically stops the refinement when a user-defined tolerance is satisfied.

\section{Acknowledgements}
Financial support provided by the Knut and Alice Wallenberg Foundation and the Swedish Research Council Grant No.
2017-04421 (VR) is gratefully acknowledged. The simulations were performed on resources provided by the Swedish National Infrastructure for Computing (SNIC) at the PDC (KTH Stockholm) and NSC (Linköping University).

\bibliographystyle{spmpsci}
\bibliography{dmas_main}

\begin{thebibliography}{10}
\providecommand{\url}[1]{{#1}}
\providecommand{\urlprefix}{URL }
\expandafter\ifx\csname urlstyle\endcsname\relax
  \providecommand{\doi}[1]{DOI~\discretionary{}{}{}#1}\else
  \providecommand{\doi}{DOI~\discretionary{}{}{}\begingroup
  \urlstyle{rm}\Url}\fi

\bibitem{akerviketal2006}
{\AA}kervik, E., Brandt, L., Henningson, D.S., Hœpffner, J., Marxen, O.,
  Schlatter, P.: Steady solutions of the {N}avier--{S}tokes equations by
  selective frequency damping.
\newblock Phys. Fluids \textbf{18}, 1--4 (2006)

\bibitem{bagheri_matrix}
Bagheri, S., {\AA}kervik, E., Brandt, L., Henningson, D.S.: Matrix-free methods
  for the stability and control of boundary layers.
\newblock AIAA J. \textbf{47}(5), 1057--1068 (2009)

\bibitem{p4est}
Burstedde, C., Wilcox, L.C., Ghattas, O.: p4est: scalable algorithms for
  parallel adaptive mesh refinement on forests and octrees.
\newblock SIAM J. Sci. Comput. \textbf{33}, 1103–1133 (2011)

\bibitem{nek5000-web-page}
Fischer, P., Lottes, J., Kerkemeier, S.: {Nek5000}.
\newblock \url{http://nek5000.mcs.anl.gov} (2008)

\bibitem{giannetti2007}
Giannetti, F., Luchini, P.: Structural sensitivity of the first instability of
  the cylinder wake.
\newblock J. Fluid Mech. \textbf{581}, 167--197 (2007)

\bibitem{huet2023}
Huet, D., Wachs, A.: A {C}artesian-octree adaptive front-tracking solver for
  immersed biological capsules in large complex domains.
\newblock J. Comput. Phys. \textbf{492}, 112424 (2023)

\bibitem{parmetis}
Karypis, G., Schloegel, K., Kumar, V.: Parmetis: Parallel graph partitioning
  and sparse matrix ordering library.
\newblock Comp. Sci. and Eng.  (1997)

\bibitem{kruse1997}
Kruse, G.W.: Parallel nonconforming spectral element solution of the
  incompressible {Navier--Stokes} equations in three dimensions.
\newblock Ph.D. thesis, {Brown University}, RI., USA (1997)

\bibitem{ARPACK}
Lehoucq, R.B., Sorensen, D.C., Yang, C.: ARPACK {U}sers' {G}uide: {S}olution of
  {L}arge-{S}cale {E}igenvalue {P}roblems with {I}mplicitly {R}estarted
  {A}rnoldi {M}ethods.
\newblock Society for Industrial and Applied Mathematics (1998)

\bibitem{luchini2014}
Luchini, P., Bottaro, A.: Adjoint equations in stability analysis.
\newblock Annu. Rev. Fluid Mech. \textbf{46}, 493--517 (2014)

\bibitem{malm2013}
Malm, J., Schlatter, P., Fischer, P.F., Henningson, D.S.: Stabilization of the
  spectral element method in convection dominated flows by recovery of
  skew-symmetry.
\newblock J. Sci. Comp. \textbf{57}, 254–277 (2013)

\bibitem{massaroStabBent}
Massaro, D., Lupi, V., Peplinski, A., Schlatter, P.: Global stability of
  180\textdegree-bend pipe flow with mesh adaptivity.
\newblock {(submitted)}  (2023)

\bibitem{massaro22TSFP}
Massaro, D., Peplinski, A., Schlatter, P.: Direct numerical simulation of
  turbulent flow around 3{D} stepped cylinder with adaptive mesh refinement.
\newblock In: Twelfth International Symposium on Turbulence and Shear Flow
  Phenomena (TSFP12) (2022)

\bibitem{massaro2023step5000}
Massaro, D., Peplinski, A., Schlatter, P.: Coherent structures in the turbulent
  stepped cylinder flow at ${Re_D}=5000$.
\newblock Int. J. Heat and Fluid Flow \textbf{102}, 109144 (2023)

\bibitem{massaro2023step1000}
Massaro, D., Peplinski, A., Schlatter, P.: The flow around a stepped cylinder
  with turbulent wake and stable shear layer.
\newblock {(submitted)}  (2023)

\bibitem{massaro22ico}
Massaro, D., Peplinski, A., Schlatter, P.: Interface discontinuities in
  spectral-element simulations with adaptive mesh refinement.
\newblock In: Spectral and High Order Methods for Partial Differential
  Equations ICOSAHOM 2020+1. Lecture Notes in Computational Science and
  Engineering, pp. 375--386. Springer International Publishing (2023)

\bibitem{KTHframework}
Massaro, D., Peplinski, A., Stanly, R., Mirzareza, S., Lupi, V., Mukha, T.,
  Schlatter, P.: {A comprehensive framework to enhance numerical simulations in
  the spectral-element code Nek5000}.
\newblock (submitted)  (2023)

\bibitem{mavriplis1989}
Mavriplis, C.: Nonconforming discretizations and a posteriori error estimators
  for adaptive spectral element techniques.
\newblock Ph.D. thesis, {MIT, USA} (1989)

\bibitem{azadinterp}
Noorani, A., Peplinski, A., Schlatter, P.: {Informal introduction to program
  structure of spectral interpolation in Nek5000}.
\newblock Tech. rep., KTH Royal Institute of Technology (2015)

\bibitem{offermans2019}
Offermans, N.: Aspects of adaptive mesh refinement in the spectral element
  method.
\newblock Ph.D. thesis, Royal Institute of Technology, KTH, Stockholm, Sweden
  (2019)

\bibitem{Offermans16}
Offermans, N., Marin, O., Schanen, M., Gong, J., Fischer, P., Schlatter, P.,
  Obabko, A., Peplinski, A., Hutchinson, M., Merzari, E.: On the strong scaling
  of the spectral element solver nek5000 on petascale systems.
\newblock In: Proceedings of the Exascale Applications and Software Conference
  2016, EASC '16. Association for Computing Machinery, New York, NY, USA (2016)

\bibitem{offermans2022}
Offermans, N., Massaro, D., Peplinski, A., Schlatter, P.: Error-driven adaptive
  mesh refinement for unsteady turbulent flows in spectral-element simulations.
\newblock Comp. Fluids \textbf{251}, 105736 (2023)

\bibitem{offermans2020}
Offermans, N., Peplinski, A., Marin, O., Schlatter, P.: Adaptive mesh
  refinement for steady flows in {N}ek5000.
\newblock Comp. Fluids \textbf{197}, 104352 (2020)

\bibitem{patera1984468}
Patera, A.T.: A spectral element method for fluid dynamics: laminar flow in a
  channel expansion.
\newblock J. Comput. Physics \textbf{54}(3), 468--488 (1984)

\bibitem{peplinski2014}
Peplinski, A., Schlatter, P., Fischer, P., Henningson, D.: Stability tools for
  the spectral-element code nek5000: Application to jet-in-crossflow.
\newblock In: Spectral and High Order Methods for Partial Differential
  Equations - ICOSAHOM 2012, pp. 349--359. Springer International Publishing
  (2014)

\bibitem{peplinski_2015}
Peplinski, A., Schlatter, P., Henningson, D.S.: Global stability and optimal
  perturbation for a jet in cross-flow.
\newblock Eur. J. Mech. B/Fluids \textbf{49}(Part B), 438--447 (2015)

\bibitem{robinson_2020}
Robinson, J.: An Introduction to Functional Analysis.
\newblock Cambridge University Press (2020)

\bibitem{schlatter_2011}
Schlatter, P., Bagheri, S., Henningson, D.S.: Self-sustained global
  oscillations in a jet in crossflow.
\newblock Theor. Comput. Fluid Dyn. \textbf{25}(1), 129--146 (2011)

\bibitem{Stability}
Schmid, P.J., Henningson, D.S.: Stability and {T}ransition in {S}hear {F}lows.
\newblock Springer (2001)

\bibitem{Sorensen_1992}
Sorensen, D.C.: Implicit application of polynomial filters in a $k$-step
  {A}rnoldi method.
\newblock SIAM J. Matrix Anal. Appl. \textbf{13}(1), 357--385 (1992)

\bibitem{Strogatz}
Strogatz, S.: Nonlinear Dynamics and Chaos: with Applications to Physics,
  Biology, Chemistry, and Engineering.
\newblock Perseus Books Group (2018)

\bibitem{Theofilis}
Theofilis, V.: Global linear instability.
\newblock Annu. Rev. Fluid Mech. \textbf{43}, 319--352 (2011)

\bibitem{Tritton}
Tritton, D.J.: Experiments on the flow past a circular cylinder at low reynolds
  numbers.
\newblock J. Fluid Mech. \textbf{6}, 547--567 (1959)

\bibitem{tuckerman2000}
Tuckerman, L.S., Barkley, D.: Bifurcation analysis for timesteppers.
\newblock In: E.~Doedel, L.S. Tuckerman (eds.) Numerical Methods for
  Bifurcation Problems and Large-Scale Dynamical Systems, pp. 453--466.
  Springer New York (2000)

\bibitem{Williamson96}
Williamson, C.H.K.: Vortex dynamics in the cylinder wake.
\newblock Annu. Rev. Fluid Mech. \textbf{28}, 477–539 (1996)

\end{thebibliography}
\end{document}